\documentclass [11] {article}

\usepackage{amsmath,amssymb,latexsym}
\usepackage[cp1251]{inputenc}
\usepackage[russian]{babel}

\newtheorem{thm}{Теорема}[section]

\def\be#1\ee{\begin{equation}#1\end{equation}}
\newcommand{\bea}{\begin{eqnarray}}
\newcommand{\eea}{\end{eqnarray}}
\newcommand{\beaa}{\begin{eqnarray*}}
\newcommand{\eeaa}{\end{eqnarray*}}

\newcommand{\bei}{\begin{itemize}}
\newcommand{\eei}{\end{itemize}}

\newcommand{\al}{\alpha}
\newcommand{\e}{\varepsilon}
\newcommand{\kk}{K}
\newcommand{\p}{\varphi}
\newcommand{\s}{\sigma}
\newcommand{\ta}{\theta}
\newcommand{\z}{\zeta}
\newcommand{\la}{\lambda}
\newcommand{\La}{\Lambda}
\newcommand{\Z}{\mathbb{Z}}
\newcommand{\X}{\mathbb{X}}
\newcommand{\PP}{\mathbb{P}}
\newcommand{\R}{\mathbb{R}}
\newcommand{\N}{\mathbb{N}}
\newcommand{\E}{\mathbb{E}}
\newcommand{\EE}{\mathcal{E}}
\newcommand{\I}{\mathbb{I}}

\begin{document}

\title{Об асимптотическом поведении\\ сложности аппроксимации
случайных полей,\\ зависящих от большого числа
параметров\thanks{Работа выполнена при финансовой поддержке РФФИ
(грант № 05-01-00911) и РФФИ-ННИО (грант № 04-01-04000). } }
\author{ Сердюкова Н.А.\footnote{Санкт-Петербургский государственный
университет, Математико-механический факультет, Библиотечная пл.,
2, 198504 Старый Петергоф, Россия.\newline
e-mail:~nora.serdyukova@ns14797.spb.edu } }

\date{}
\maketitle
\bigskip

\begin{abstract}
В настоящей статье изучено поведение сложности аппроксимации в
среднем для $d$-параметрических случайных полей тензорного типа. В
работе \cite{LT} было показано, что для заданного уровня
относительной ошибки сложность аппроксимации возрастает
экспоненциально при $d\to \infty$, то есть наблюдается так
называемый феномен "проклятия размерности".
\par
В данной статье для сложности аппроксимации получено точное
асимптотическое выражение.
\end{abstract}

\noindent

{\bf Ключевые слова:}\ случайные поля, гауссовские процессы,
ошибка линейной аппроксимации, сложность, проклятие размерности.
\bigskip\bigskip
\baselineskip=6.0mm

\section{Введение}

\par Пусть $T$ -- некоторое параметрическое множество.
Предположим, что случайная функция $ X(t), \; t \in T $ допускает
представление в виде ряда \[ X(t) = \sum_{k=1}^{\infty}\xi_k
\p_k(t),\] где $\xi_k$ -- случайные величины, а $\p_k$ --
детерминированные вещественные функции.

Для каждого конечного множества натуральных чисел $\kk\subset \N$
обозначим $ X_\kk(t) = \sum_{k\in \kk} \xi_k \p_k(t)$. При решении
многих задач бывает необходимо аппроксимировать $X$ процессом
конечного ранга $ X_\kk$, в связи с чем возникают естественные
вопросы: какого размера должно быть выбрано множество $\kk$, чтобы
обеспечить заданную точность аппроксимации? Среди всех множеств
$\kk$ заданного размера, множество какого вида обеспечит
наименьшую ошибку?

\par В настоящей статье мы
ответим на первый из поставленных вопросов применительно к
определенному классу случайных функций,  а именно, {\it случайным
полям тензорного типа} с параметрическими множествами высокой
размерности.

\medskip
Пусть неотрицательная последовательность $(\lambda(i))_{i\ge 1}$
удовлетворяет условию \be  \label{l2}
  \sum_{i=1}^\infty \lambda(i)^2 <\infty
\ee и пусть функции $(\varphi_i)_{i>0}$ образуют полную ортонормированную
систему в $L_2[0,1]$.

\par Рассмотрим семейство случайных полей тензорного типа
\[\X= \left\{ X^{(d)}(t), t \in[0,1]^d \right\},\; d=1,2,...,\]
заданных формулой

\begin{equation} \label{xd1}
X^{(d)}(t)= \sum_{k\in \N^d} \xi_k \ \prod_{l=1}^d \lambda(k_l)
\prod_{l=1}^d \varphi_{k_l}(t_l),
\end{equation}

\noindent где $\xi_k$ -- некоррелированные случайные величины с
нулевым средним и единичной дисперсией. Очевидно, что если
выполнено условие \eqref{l2}, то траектории $X^{(d)}$ принадлежат
$L_2([0,1]^d)$ почти наверное, и система собственных чисел
ковариационного оператора $X^{(d)}$ имеет вид
\begin{equation}\label{sobd}
    \lambda_k^2:= \prod_{l=1}^{d} \la(k_l)^2\,,\,\,\,\, k \in \N^d.
\end{equation}
В дальнейшем мы не упоминаем индекс $d$, а просто пишем
$X(t)$ вместо $X^{(d)}(t)$.

\par Пусть $T=[0,1]^d$. Для любого
$n>0$ обозначим через $X_n$ частичные суммы ряда (\ref{xd1}),
отвечающие $n$ максимальным собственным числам. Нас интересует
асимптотическое поведение ошибки аппроксимации в среднем $X$
посредством $X_n$, а именно \[
       e(X, X_n; d) = \left(\E ||X-X_n ||^2_{L_2(T)}
       \right)^{1/2}\;,\;\;\;\;\;d \to \infty.
\] Далее мы будем
рассматривать только $L_2(T)$-нормы, поэтому под символом $||
\cdot||$ всегда будем подразумевать $||\cdot ||_{L_2(T)}$. Хорошо
известно (см., например, \cite{BS}, \cite{KL} или \cite{R}), что
среди всех линейных аппроксимаций порядка $n$, $X_n$ обеспечивает
минимальную среднеквадратичную ошибку.

\par Поскольку мы собираемся исследовать {\it семейство} случайных функций,
то более естественно изучать {\it относительные}
ошибки, то есть сравнивать величину ошибки с величиной самой
функции.

\par Обозначим
\[ \La:=\sum_{i=1}^\infty \la(i)^2, \]
тогда
\[
\E\| X\|^2 = \sum_{k\in\N^d} \la_k^2 = \La^d.
\]

\par Определим относительную сложность среднеквадратичной
аппроксимации как
\[
n(\e,d) :=
 \min \{n : \frac{e(X, X_n; d)}{\left(\E\| X\|^2\right)^{1/2}} \leq \e\}
= \min \{n : \E \|X-X_n \|^2 \leq \e^2\La^d \}.
\]

Заметим, что интересующее нас поведение сложности аппроксимации
$n(\e,d)$, принадлежит к классу проблем, связанных с изучением
зависимости сложности линейных многопараметрических задач от
размерности, см. работы Х.~Возняковского (\cite{W92}, \cite{W94a},
\cite{W94b}, \cite{W2006}) и приведенные там ссылки.

\par Для изучения свойств массива собственных чисел \eqref{sobd}
(детерминированного!) в \cite{LT} было предложено использовать
вспомогательную вероятностную конструкцию. Мы также будем
придерживаться этого подхода.
\par Рассмотрим последовательность независимых одинаково
распределенных случайных величин $\left(U_l\right),
\,\,l=1,2,...$, общее распределение которых задано формулой
\be\label{def_Ul} \PP(U_l=-\log \la(i))=\frac{\la(i)^2}{\La}
\,,\,\,\,\,i=1,2,... \ee
\par Если  выполнено условие \be\label{3d_moment} \sum_{i=1}^{\infty}|\log
\la(i)|^3 \la(i)^2 \;<\; \infty \ee то, очевидно,  $ \E | U_l |^3
< \infty.$

\par Обозначим через $M$ и $\s^2$ математическое ожидание и дисперсию
$U_l$, соответственно. Ясно, что  \bea \nonumber M &=& -
\sum_{i=1}^{\infty}\log \la(i)\,\frac{ \la(i)^2}{\La},
\\ \nonumber
\s^2 &=& \sum_{i=1}^{\infty}|\log \la(i)|^2\, \frac{
\la(i)^2}{\La}\; -\; M^2. \eea Тогда третий центральный момент
$U_l$ задается
\[
\al^3 := \E(U_l -M)^3 =  - \sum_{i=1}^{\infty}\log
\la(i)^3\,\frac{ \la(i)^2}{\La} \; - \; 3M\s^2 \;-\;M^3.
\]
Если выполнено \eqref{3d_moment}, то $|M|<\infty$, $0 \leq \s^2 <
\infty$ и $|\al|< \infty.$

\par В дальнейшем "взрывающийся" коэффициент \be\label{Expl} \EE :=
\La e^{2M} \ee будет играть значительную роль. В~\cite{LT}
показано, что $\EE >1$ за исключением вырожденного случая, когда
число положительных собственных чисел равно нулю или единице.
Другими словами, $\EE =1$ тогда и только тогда, когда $\s=0$. Мы
исключим этот случай из дальнейшего рассмотрения.

В \cite{LT} получен следующий результат (теорема 3.2).

\begin{thm}\label{T:LT_3.2}
Предположим, что последовательность $(\la(i))$, $i=1,2,...$
удовлетворяет условию
\[  \sum_{i=1}^{\infty}|\log \la(i)|^2\,
\la(i)^2 < \infty.
\]
Тогда для каждого $\e \in (0,1)$ выполнено
\[
\lim_{d \to \infty} \frac{\log n(\e,d) - d \log \EE}{\sqrt{d}} =
2q,
\]
где квантиль $q=q(\e)$ выбрана из уравнения \be\label{q} 1 -
\Phi\left( \frac{q}{\s}\right) = \e^2. \ee
\end{thm}

Авторы \cite{LT} предположили, что при более сильных условиях на
последовательность $(\la(i))$ можно доказать, что
\[
n(\e,d) \approx \frac{C(\e) \EE^d e^{2q\sqrt{d}}}{\sqrt{d}}\ ,\;
\;\, d \to \infty.
\]

Мы покажем, что выполнено даже несколько более точное утверждение.

\section{Основной результат}

\par Нам придется отдельно рассматривать два случая
в зависимости от природы распределения величин $U_l$. Доказательство и
конечный результат будут зависеть от того, является это распределение
{\it решетчатым} или нет.
\par Напомним, что дискретное распределение случайной величины $U$ называется
решетчатым, если существуют такие числа $a$ и $h>0$, что все
возможные значения $U$ могут быть представлены в виде $a + \nu h$,
где $\nu$ принимает целые значения. Число $h$ называется шагом
распределения. В дальнейшем при рассмотрении решетчатого случая
мы предполагаем, что $h$ является максимальным шагом
распределения, то есть что нельзя представить все возможные
значения $U_l$ в виде $b + \nu h_1$ с некоторыми $b$ и $h_1>h$.

\par Из определения \eqref{def_Ul} следует, что величины $U_l$ имеют
решетчатое распределение тогда и только тогда, когда
$\la(i) = C e^{-n_i h}$ для некоторых положительных вещественных
$C$, $h$ и $n_i \in \N$. Мы будем называть эту ситуацию {\it
решетчатым случаем} и будем предполагать, что $h$ выбрано
максимально возможным. В противном случае мы будем говорить, что
имеет место {\it нерешетчатый случай}.

\bigskip
\begin{thm}\label{T}
Пусть последовательность $\left( \la(i)\right)$, $i=1,2,...$
удовлетворяет условию~\eqref{3d_moment}. \par Тогда для каждого
$\e \in (0,1)$ выполнено
\[ n(\e,d) = K\ \phi(\frac{q}{\s})\  \EE^d e^{2q \sqrt{d}}
\,d^{-1/2} \left( 1 + o(1) \right),\; \;\, d \to \infty,\] где
\beaa\phi(x) &=& \frac{1}{\sqrt{2 \pi}} \, e^{-x^2/2 },\\
K &=& \begin{cases}
\frac{h}{\s(1-e^{-2h})} & \text{в решетчатом случае;} \\
 \frac{1}{2\s}  & \text{в нерешетчатом;}
\end{cases}
\eeaa и квантиль $q=q(\e)$ определена в~\eqref{q}.
\end{thm}

\bigskip

\begin{remark}
\bei \item Мы видим, что сложность аппроксимации возрастает
экспоненциально когда $d \to \infty$. Это явление обычно называют
"проклятием размерности" (dimensionality curse) или
"intractability".  См., например,~\cite{R} и~\cite{W94a}. Понятие
"проклятия размерности"  восходит, по крайнай мере, к работам
Беллмана~\cite{B}. \item По правилу Лопиталя \[ \lim_{h \to 0}
\frac{h}{1-e^{-2h}} = \frac{1}{2 }\ , \] так что формулы для $K$
согласованы между собой при $h\to 0$. \eei
\bigskip
\end{remark}

\begin{pf}

Определим $\z = \z(\e, d)$ как максимальное положительное
вещественное число такое, что сумма собственных чисел, меньших чем
$\z^2$, удовлетворяет неравенству
\[ \sum_{k\in\N^d : \la_k< \z}\la_k^2
\leq \e^2\La^d.
\]
Определим решетчатое множество в  $\N^d$
\[
\mathrm{A} =  \mathrm{A}(\e, d):= \{  k \in \N^d : \la_k \geq \z
\} = \{  k \in \N^d : \prod_{l=1}^{d} \la(k_l) \geq \z \}.
\]
Поскольку для любого $k \in \mathrm{A}$ верно, что $ \la_k
>0$, то  \beaa \nonumber && n(\e,d) \; = \; card
\left(\mathrm{A}\right) \;
 =\;
\sum_{k \in \mathrm{A}}\frac{\la_k^2}{\la_k^2}
\\ \nonumber
&=& \sum_{ k \in \N^d : -\sum \log \la(k_l) \leq -\log \z}\La^d
\exp\{-2 \sum_{l=1}^d \log \la(k_l)\}\prod_{l=1}^{d}\PP(U_l=-\log
\la(k_l))
\\ \nonumber
&=& \La^d \,\E\exp \{2 \sum_{l=1}^d U_l  \} \I_{\{\sum_{l=1}^d U_l
\leq -\log \z \}}. \eeaa Для центрированных и нормированных сумм
\[ Z_d = \frac{\sum_{l=1}^d U_l
- d M }{\s\sqrt{d}}
\]
выполнено
\[ \{\sum_{l=1}^d U_l \leq -\log
\z \} = \{ Z_d \leq \theta \},
\]
где \be\label{ta_def} \ta = \ta(\e,d) = - \frac{\log\z + d M}{\s
\sqrt{d}}. \ee

Покажем теперь, что $\ta$ допускает полезную вероятностную
интерпретацию в терминах случайных величин $U_l$ и их сумм.
Действительно, применяя лемму~3.1 из~\cite{LT}, мы имеем для любых
$d \in \N$ и $z \in \R^1$
 \beaa
 \sum_{k\in \N^d :\la_k < z} \la_k ^2 \; & = & \La^d\;\PP \left( \sum_{l=1}^d U_l > -\log z
 \right)\\ =  \La^d\;\PP \left( Z_d > - \frac {\log z + d M}{\s
 \sqrt{d}}\right) &=&  \La^d\;\PP \left( Z_d > \ta_z   \right),
 \eeaa
 где
\[ \ta_z = -\frac{\log z + d M}{\s \sqrt{d}}.
\]
Зафиксируем $\e \in (0,1)$. Заметим, что \[   \sum_{k\in \N^d
:\la_k < z} \la_k ^2 \; \leq \; \e^2 \La^d
\]
тогда и только тогда, когда \[  \PP \left( Z_d > \ta_z   \right)
\; \leq \; \e^2 .
\]
Поэтому $\ta =\ta(\e, d)$, определенная равенством~\eqref{ta_def},
является $(1-\e^2)$-квантилью распределения $Z_d$, а именно,
\[
\ta(\e, d)=\min\{\ta:\ \PP\left(Z_d>\ta\right)\le \e^2\}
 =\min\{\ta:\ \PP\left(Z_d\le \ta\right) >1- \e^2\}.
\]
Обозначим через $q=q(\e)$ квантиль функции распределения
нормального закона, выбранную из уравнения~\eqref{q}. Тогда, в
силу центральной предельной теоремы, \be\label{ta}
\ta(d,\e) \to \frac{q(\e)}{\s}\; , \;\; d \to \infty, \ee для
каждого фиксированного $\e \in(0,1)$.

Теперь мы можем вернуться к изучению сложности аппроксимации. Мы
получили, что
 \bea \nonumber n(\e,d) &=& \EE^d \, \E \exp
\{ 2 \s \sqrt{d} Z_d\} \I_{\{Z_d \leq \ta \}}
\\ \nonumber
&=& \EE^d \, \exp \{ 2 \s \sqrt{d}\ta \} \int_{-\infty}^{\ta} \exp
\{ 2 \s \sqrt{d}(z-\ta) \} \,\mathrm{d}F_d(z) , \eea где $F_d(z) =
\PP(Z_d < z)$ и $\EE$ определена в~\eqref{Expl}.

Обозначим
\[
  \Psi_d(z):= \exp \{ 2 \s \sqrt{d}(z-\ta) \}
\]
и проинтегрируем по частям интеграл
\[
 \int_{-\infty}^{\ta} \Psi_d(z)  \,\mathrm{d}[F_d(z)-F_d(\ta)] =
 \int_{-\infty}^{\ta} [- F_d(z)+F_d(\ta)] \,\mathrm{d}\Psi_d(z).
\]
\vspace{0.3cm}

\begin{bf}  Нерешетчатый случай \end{bf}

В этой части доказательства мы будем предполагать, что
распределение $\left( U_l \right)$ не является решетчатым. Это
имеет место в наиболее интересных случаях, таких как броуновский
лист (поле Винера-Ченцова), дробный броуновский лист, броуновская
подушка, поле Кифера.

\par В силу того, что выполнено~\eqref{3d_moment} мы можем применить
теорему Крамера-Эссеена (см. теорему~2 \S 42 в~\cite{GK},
теорему~21 \S 7 главы~V в~\cite{Pe2} либо теорему~4 \S 3  главы~VI
в~\cite{Pe1}), откуда немедленно получаем \bea \label{CE_th}
\nonumber \int_{-\infty}^{\ta} [- F_d(z)+F_d(\ta)]
\,\mathrm{d}\Psi_d(z)& =& \int_{-\infty}^{\ta} [-
\Phi(z)+\Phi(\ta)] \,\mathrm{d}\Psi_d(z)
\\ \nonumber
+ \frac{\al^3}{6\s^3 \sqrt{2\pi d}}\int_{-\infty}^{\ta}
[(z^2-1)e^{-z^2/2} &-& ((\ta^2-1)e^{-\ta^2/2} ]
\,\mathrm{d}\Psi_d(z) + o\left( \frac{1}{\sqrt{d}}\right)
\\
= I_1 + I_2 -I_3  - I_4 &+& o\left( \frac{1}{\sqrt{d}}\right),
\eea где \beaa \nonumber I_1 &=& \int_{-\infty}^{\ta} [-
\Phi(z)+\Phi(\ta)] \,\mathrm{d}\Psi_d(z),
\\
I_2 &=& \frac{\al^3}{6\s^3 \sqrt{2\pi d}}\int_{-\infty}^{\ta}  z^2
e^{-z^2/2}  \,\mathrm{d}\Psi_d(z),
\\ \nonumber
I_3 &=& \frac{\al^3}{6\s^3 \sqrt{2\pi d}}\int_{-\infty}^{\ta}
e^{-z^2/2} \,\mathrm{d}\Psi_d(z),
\\ \nonumber
I_4 &=& \frac{\al^3}{6\s^3 \sqrt{2\pi d}}\left(\ta^2 - 1 \right)
e^{-\ta^2/2} \sim \frac{\al^3}{6\s^3 \sqrt{2\pi
d}}\left(\left(\frac{q}{\s}\right)^2 - 1 \right)
e^{-q^2/2\s^2}\;,\; d\to \infty. \eeaa Последнее соотношение
следует из~\eqref{ta}. Поскольку $\mathrm{d}\Psi_d(z) = 2 \s
\sqrt{d}\Psi_d(z)d z$,  то  после замены переменных $I_2$ примет
вид \beaa \nonumber
 I_2 & = & I_2(d,\ta) =
 \frac{\al^3}{3\s^2 \sqrt{2
 \pi d }}\int_0^{\infty}(\ta - \frac{y}{\sqrt{d}})^2 \exp\{
 -\frac{1}{2}(\ta - \frac{y}{\sqrt{d}})^2\}\exp\{-2\s y \} \; \mathrm{d} y
 \eeaa
 где $y = -\sqrt{d}(z-\ta)$.

\par  Для любого $ d =1,2,...$ верно неравенство

\[
0\le \left(\ta - \frac{y}{ \sqrt{d}}\right)^2 \exp\{
-\frac{1}{2}(\ta - \frac{y}{ \sqrt{d}})^2 \} \leq (|\ta| + y)^2.
\]
Эта оценка дает нам суммируемую мажоранту, необходимую для
применения теоремы Лебега. Учитывая \eqref{ta} и переходя к
пределу под знаком интеграла, мы получаем при $d \to \infty$,
\[
I_2(d,\ta) =
 \frac{\al^3}{6\s^3 \sqrt{2
 \pi d }} \;\left(\frac{q}{\s}\right)^2 e^{-q^2/2\s^2} \left( 1 +
 o(1)\right).
\]
Аналогично,
\[
 I_3(d,\ta) =
 \frac{\al^3}{6\s^3 \sqrt{2
 \pi d }}  e^{-q^2/2\s^2} \left( 1 + o(1)\right).
\]
Таким образом, мы получили, что $ \sqrt{d}I_4 =  \sqrt{d}(I_2 -I_3
)\left( 1 + o(1)\right)$, следовательно, $I_2-I_3-I_4 =  o\left(
\frac{1}{\sqrt{d}}\right)$.

Рассмотрим основной интеграл $I_1$.

\bea \label{I1} \nonumber I_1 &=& I_1(d, \ta) = \int_{-\infty}^\ta
[-\Phi(z)+\Phi(\ta)] \,\mathrm{d}\Psi_d(z)\\
\nonumber &=&\frac{1}{\sqrt{2 \pi}} \int_{-\infty}^\ta  \exp\{ 2 \s \sqrt{d}(z-\ta) \} \exp\{ -z^2/2\}  \; \mathrm{d} z \\
\nonumber &=&\frac{1}{\sqrt{2 \pi d }}\int_0^{\infty}\exp\{
 -\frac{1}{2}(\ta - \frac{y}{\sqrt{d}})^2\}\exp\{-2\s y \} \; \mathrm{d} y
\\ &\sim&  \frac{1}{2\s \sqrt{2 \pi d}} e^{-q^2/2\s^2}\;,\;\; d \to \infty. \eea

 Таким образом, мы получили желаемую оценку
\[
n(\e,d) = \frac{\EE^d\, \exp\{ 2q \sqrt{d} \}}{2 \s \sqrt{d}}
\frac{1}{\sqrt{2 \pi}} \exp\{-q^2/2\s^2 \} \left( 1 + o(1)
\right).
\]
\vspace{0.3cm}

\begin{bf} Решетчатый случай \end{bf}

Теперь мы будем действовать в предположении, что
последовательность $\left(U_l\right)$ имеет решетчатое
распределение.

\par Пусть случайная величина $U_l$ принимает значения следующего
вида
\[
\tilde{a} + \nu h , \; \nu = 0,\pm 1, \pm 2,...
\] где
$\tilde{a} = M+a$ -- сдвиг, а  $h$ -- максимальный шаг
распределения. Тогда все возможные значения суммы $Z_d$
представимы в виде
\[
\frac{d a + \nu h}{\s \sqrt{d}}, \; \nu = 0,\pm 1, \pm 2,...
\]
Введем в рассмотрение функцию \[ S(x) = [x] - x+ \frac{1}{2}, \]
где $[x]$ обозначает, как обычно, целую часть $x$, и определим
\[
 S_d(x) = \frac{h}{\s} \,S\left( \frac{x \s \sqrt{d} - d
a}{h}\right).
\]
Пусть $F_d(z)$ определена как прежде. Тогда, если выполнено
\eqref{3d_moment}, то результат Эссеена (см. теорему~1~\S~43
в~\cite{GK}) влечет равенство
\[
F_d(z) - \Phi(z) = \frac{e^{-z^2/2}}{\sqrt{2\pi}} \left(
\frac{S_d(z)}{\sqrt{d}} - \frac{\al^3 (z^2 - 1)}{6 \s^3
\sqrt{d}}\right) + o \left( \frac{1}{\sqrt{d}}\right)
\]
равномерно по $z$.

Сравнивая с \eqref{CE_th}, мы видим, что  нам нужно оценить только
дополнительное слагаемое \beaa J &=&\frac{1}{\sqrt{2 \pi d}}
\int_{-\infty}^{\ta}[- S_d(z) e^{-z^2/2} + S_d(\ta) e^{-\ta^2/2}]
\mathrm{d} \Psi_d(z)
\\
&=& \frac{1}{\sqrt{2 \pi d}} \int_{-\infty}^{\ta} \Psi_d(z)
\mathrm{d} \left( S_d(z) e^{-z^2/2}\right) =  J_1 - J_2 + J_3 ,
\eeaa где \beaa \nonumber J_1 &=& \frac{1}{\sqrt{2 \pi d}}
\int_{-\infty}^{\ta} \Psi_d(z) S_d'(z)e^{-z^2/2} \mathrm{d} z,
\\ \nonumber J_2 &=& \frac{1}{\sqrt{2 \pi d}}
\int_{-\infty}^{\ta} \Psi_d(z) S_d(z) z e^{-z^2/2} \mathrm{d} z,
\eeaa  и $J_3$ -- "дискретная часть", которая определена следующим
образом.

Заметим, что $S(x)$ является периодической функцией с единичным
периодом, поэтому $S_d(x)$ обладает периодом $h/\s \sqrt{d}$ и
имеет скачки в точках вида $\{ \frac{kh + da}{\s \sqrt{d}}, k \in
\Z \}$. Если $\ta$ принадлежит этой решетке, то существует целое
$k'$ такое, что $\ta = \frac{k'h + da}{\s \sqrt{d}}$.
Следовательно, можно интегрировать разрывную часть интеграла $J$
по отношению к дираковской мере $\frac{h}{\s}\delta_{\frac{kh +
da}{\s \sqrt{d}}}$. Тогда \[J_3 = \frac{1}{\sqrt{2 \pi d}}\,
\frac{h}{\s} \sum_{k=-\infty}^{k'} \Psi_d\left(\frac{kh + da}{\s
\sqrt{d}}\right) \exp\{-\frac{1}{2}\left(\frac{kh + da}{\s
\sqrt{d}}\right)^2\}.\]

\par Оценим сперва интеграл $J_1$. В тех точках, где производная
$S_d'(z)$ имеет смысл, можно легко вычислить, что $S_d'(z) =
\frac{h}{\s}S\left(\frac{z\s\sqrt{d} - da}{h} \right) = -\sqrt{d}
$, тогда, подобно нерешетчатому случаю, по теореме Лебега мы имеем
\bea \label{J1}\nonumber J_1 &=& \frac{-\sqrt{d}}{\sqrt{2 \pi d}}
\int_{-\infty}^{\ta} \exp\{ 2\s \sqrt{d}(z-\ta)\} \exp\{ -z^2/2\}
\mathrm{d} z
\\\nonumber &=& \frac{-1}{\sqrt{2 \pi d}} \int_{0}^{\infty} \exp
\{-\frac{1}{2}(\ta
- \frac{y}{\sqrt{d}})^2 \}\exp\{-2\s y \} \mathrm{d} y\\
&\sim& \frac{-1}{ 2\s \sqrt{2 \pi d} }e^{-q^2/2\s^2}\;,\; d \to
\infty , \eea следовательно, $\sqrt{d}J_1 = - \sqrt{d} I_1\left(
1+ o(1)\right)$.

Что касается интеграла $J_2$, то он, при достаточно больших $d$,
вообще не играет роли. Действительно, \beaa J_2 &=&
\frac{1}{\sqrt{2 \pi d}} \int_{-\infty}^{\ta}
\exp\{ 2\s \sqrt{d}(z-\ta)\}S_d(z) z  \exp\{ -z^2/2\} \mathrm{d} z \\
 &=&       \frac{1}{\sqrt{2 \pi d}}  \frac{1}{\sqrt{d}} \int_{0}^{\infty} \exp \{-\frac{1}{2}(\ta
- \frac{y}{\sqrt{d}})^2 \} (\ta - \frac{y}{\sqrt{d}})S_d(\ta
- \frac{y}{\sqrt{d}})\exp\{-2\s y \} \mathrm{d}y \\
&\leq& \frac{3 h}{2 \s d\sqrt{2 \pi } }\int_{0}^{\infty} \exp
\{-\frac{1}{2}(\ta - \frac{y}{\sqrt{d}})^2 \} (\ta -
\frac{y}{\sqrt{d}})\exp\{-2\s y \} \mathrm{d}y
\\ &\sim& \frac{3h}{ 4\s^2 d \sqrt{2 \pi }  }\left( \frac{q}{\s}\right)^2 e^{-q^2/2\s^2}\;,\; d \to \infty
... \eeaa И, разумеется,  $J_2 = o\left(
\frac{1}{\sqrt{d}}\right)$.

\par Теперь мы перейдем к рассмотрению наиболее важного слагаемого
 \bea \label{J3} \nonumber J_3 &=&
\frac{1}{\sqrt{2 \pi d}}\, \frac{h}{\s} \sum_{k=-\infty}^{k'}
\exp\{2\s \sqrt{d} \left( \frac{k h + d a }{\s \sqrt{d}} - \ta
\right) \} \exp\{- \frac{1}{2} \left(
 \frac{k h + d a}{\s \sqrt{d}} \right)^2\}
\\\nonumber&=& \frac{1}{\sqrt{2 \pi d}}\,
\frac{h}{\s} \sum_{k=-\infty}^{k'}\exp\{2 h (k - k') \} \exp\{-
\frac{1}{2} \left(
 \frac{k h + d a}{\s \sqrt{d}} \right)^2\}
 \\ \nonumber&=& \frac{1}{\sqrt{2 \pi d}}\,
\frac{h}{\s} \sum_{l=0}^{\infty}  \exp\{-2 h l\} \exp\{-
\frac{1}{2} \left(
 \frac{(k' - l) h + d a}{\s \sqrt{d}} \right)^2\}
 \\\nonumber &=& \frac{1}{\sqrt{2 \pi d}}\,
\frac{h}{\s} \sum_{l=0}^{\infty}  \exp\{-2 h l\} \exp\{-
\frac{1}{2} \left( \ta - \frac{l h}{\s \sqrt{d}}\right)^2
\\ & \sim& \frac{1}{\s
\sqrt{ d}}\, \frac{h }{(1-e^{-2h})} \frac{1}{\sqrt{2
\pi}}e^{-q^2/2\s^2}\;,\; d \to \infty. \eea Таким образом, мы
получили \[ \sqrt{d}J_3 = \sqrt{d} \frac{2h }{(1-e^{-2h})} I_1
\left( 1 + o(1)\right).\]

\par Объединяя вместе оценки \eqref{I1}, \eqref{J1} и \eqref{J3},
мы получаем  \[ n(\e,d) = \frac{\EE^d \, e^{ 2q \sqrt{d} }}{ \s
\sqrt{d}} \frac{h }{(1-e^{-2h})} \frac{1}{\sqrt{2 \pi}} \,
e^{-q^2/2\s^2 } \left( 1 + o(1) \right),\; \;\, d \to \infty. \]

\end{pf}
Теорема~\ref{T} доказана.

\bigskip

\noindent Настоящая статья была частично написана во время
пребывания автора в Институте математической стохастики
университета Георга-Августа, Г\"{е}ттинген. Хочется особенно
поблагодарить профессора М.А. Лифшица за постановку данной задачи
и постоянную поддержку, а также профессора М. Денкера за
содействие и обеспечение прекрасных условий для работы.

\bibliographystyle{amsplain}

\bigskip
\bigskip
{\footnotesize \baselineskip=12pt \noindent
\parbox[t]{5cm}

\end{document}